\documentclass[psamsfonts,reqno]{amsart}
\usepackage{amssymb}
\usepackage{amscd}
\usepackage{latexsym}

\theoremstyle{plain}
\newtheorem{lemma}{Lemma}[section]
\newtheorem{theorem}[lemma]{Theorem}
\newtheorem{proposition}[lemma]{Proposition}
\newtheorem{corollary}[lemma]{Corollary}

\newtheorem*{sclaim}{Claim}
\newtheorem*{thm}{Theorem 1}
\newtheorem*{tha}{Theorem 2}
\newtheorem*{thb}{Theorem 3}

\theoremstyle{definition}
\newtheorem{definition}[lemma]{Definition}

\newtheorem{example}[lemma]{Example}
\newtheorem*{problem}{Problem}
\newtheorem*{remark}{Remark}

\numberwithin{equation}{section}

\newcommand{\qedc}{{\qed}~{\rm Claim~{\theclaim}.}}
\newcommand{\sqedc}{{\qed}~{\rm Claim.}}

\newenvironment{scproof} {\begin{proof}[Proof of Claim.]}
{\sqedc\renewcommand{\qed}{}\end{proof}}

\DeclareMathOperator{\Conv}{Conv}
\DeclareMathOperator{\Min}{Min}

\newcommand{\upw}{\mathbin{\uparrow}}

\newcommand{\cm}{commutative monoid}

\newcommand{\poag}{partially ordered abelian group}
\newcommand{\povs}{partially ordered right vector space}
\newcommand{\nor}{normally ordered ring}

\newcommand{\xx}{\mathbf{x}}

\newcommand{\NN}{\mathbb{N}}
\newcommand{\ZZ}{\mathbb{Z}}
\newcommand{\QQ}{\mathbb{Q}}
\newcommand{\RR}{\mathbb{R}}

\newcommand{\eC}{\ensuremath{\mathbf{C}}}
\newcommand{\eP}{\ensuremath{\mathbf{P}}}
\newcommand{\eI}{\ensuremath{\mathbf{I}}}
\newcommand{\eE}{\ensuremath{\mathbf{E}}}

\renewcommand{\SS}{\mathcal{S}}

\def\vv<#1>{\langle#1\rangle}
\def\fs(#1){\{1,\ldots,#1\}}

\begin{document}

\title[Finitely presented ordered groups]
{Finitely presented, coherent, and ultrasimplicial\\
ordered abelian groups}

\author{J.~F. Caillot}
\address{D\'epartement de Math\'ematiques\\
Universit\'e de Caen\\
14032 Caen Cedex\\
France}
\email{caillot@math.unicaen.fr}

\author{F. Wehrung}
\address{C.N.R.S., E.S.A. 6081\\
D\'epartement de Math\'ematiques\\
Universit\'e de Caen\\
14032 Caen Cedex\\
France}
\email{gremlin@math.unicaen.fr}
\urladdr{http://www.math.unicaen.fr/\~{}wehrung}


\begin{abstract}
We study notions such as finite presentability and coherence, for
partially ordered abelian groups and vector spaces. Typical results
are the following:
\begin{enumerate}
\item A \poag\ $G$ is finitely presented if and only if $G$ is
finitely generated as a group, $G^+$ is well-founded as a partially
ordered set, and the set of minimal elements of $G^+\setminus\{0\}$ is
finite.

\item Torsion-free, finitely presented \poag s can be represented as
subgroups of some $\ZZ^n$, with a finitely generated submonoid of
$(\ZZ^+)^n$ as positive cone.

\item Every unperforated, finitely presented \poag\ is Archimedean.
\end{enumerate}
Further, we establish connections with interpolation. In
particular, we prove that a divisible dimension group $G$ is a directed
union of simplicial subgroups if and only if every finite subset of $G$
is contained into a finitely presented ordered subgroup.
\end{abstract}

\maketitle

\section*{Introduction}
The elementary theory of abelian groups implies that every finitely
generated abelian group is finitely presented. The situation for
\emph{partially ordered} abelian groups is very different, as easy
examples show. In \cite{Wehr}, the second author establishes a general
framework for convenient study of notions related to finite
presentability for partially ordered modules over partially ordered
rings. This framework yields for example statements such as
the following:

\begin{thm}
Let $G$ be a \poag. Then the following statements hold:
\begin{enumerate}
\item $G$ is finitely presented if and only if $G$ is finitely
generated as a group, and $G^+$ is finitely generated as a monoid.

\item If $G$ is finitely presented, then every ordered subgroup of $G$
is finitely presented.
\end{enumerate}
\end{thm}

In fact, this result can be extended to the much more general
context of partially ordered right modules over a \emph{right coherent}
\nor. This context covers the case of \povs s over a given totally
ordered field, or, more generally, \emph{division ring}, since
commutativity is not required. In particular, the methods used are
absolutely not specific to groups.

By contrast, we obtain here specific results, such as the following
(see Theorem~\ref{T:CharFPgroup}):

\begin{tha}
Let $G$ be a \poag. Then $G$ is finitely presented, if and only if $G$
is finitely generated as a group, $G^+$ is well-founded (for its
natural ordering), and the set of minimal elements of
$G^+\setminus\{0\}$ is finite.
\end{tha}

We also obtain a general representation result of torsion-free,
finitely presented \poag s, as groups of the form $\ZZ^n$, endowed
with a finitely generated submonoid of $(\ZZ^+)^n$ as positive cone
(see Theorem~\ref{T:DescFPgp}). As a corollary, we obtain
(Corollary~\ref{C:NonPerfAr}) that \emph{every unperforated, finitely
presented \poag\ is Archimedean}.

The theory developed is convenient to study a variant of a notion of
\emph{ultrasimpliciality}, introduced by G.~A. Elliott, see
\cite{Elli1,Elli2}. Say that a \poag\ $G$ is \emph{simplicial}, if it is
isomorphic to some $\ZZ^n$ endowed with the natural ordering, and
\emph{\eE-ultrasimplicial}, if it is a directed union of simplicial
subgroups. We prove, in particular, the following result (see
Corollary~\ref{C:EUltrgp}):

\begin{thb}
A divisible dimension group is \eE-ultrasimplicial if and only if it is
\emph{coherent}, that is, every finitely generated ordered subgroup is
finitely presented.
\end{thb}

We also establish an analogue of this result for \povs s over totally
ordered division rings (see Corollary~\ref{C:EUltr}).

\section{Basic concepts}\label{S:BasConc}

Let us first recall some basic definitions, see \cite{Gpoag}.
A \emph{\poag} is, by definition, an abelian group $\vv<G,+,0>$,
endowed with a partial ordering $\leq$ that is \emph{translation
invariant}, that is, $x\leq y$ implies
$x+z\leq y+z$, for all $x$, $y$, $z\in G$. We put, then
 \[
 G^+=\{x\in G\mid x\geq0\},\quad\text{and}\quad
 G^{++}=G^+\setminus\{0\}.
 \]
If $G$ and $H$ are \poag s, a \emph{positive homomorphism} from $G$ to
$H$ is an order-preserving group homomorphism from $G$ to $H$. A
positive homomorphism $f$ is a \emph{positive embedding}, if $f(x)\geq0$
if and only if $x\geq0$, for all $x\in G$. We shall use the notation
$\NN=\ZZ^+\setminus\{0\}$.

We say that a \poag\ $G$ satisfies the
\emph{(Riesz) interpolation property} if, for all $a_0$, $a_1$, $b_0$,
and $b_1$ in $G$, if $a_0,a_1\leq b_0,b_1$, then there exists $x\in G$
such that $a_0,a_1\leq x\leq b_0,b_1$, and we say then that $G$ is an
\emph{interpolation group}. We say that $G$ is \emph{directed} if
$G=G^++(-G^+)$, \emph{unperforated} if $mx\geq0$ implies that $x\geq0$,
for all $m\in\NN$ and all $x\in G$. A
\emph{dimension group} is, by definition, a directed, unperforated
interpolation group.

An important particular class of dimension groups is the class of
\emph{simplicial groups}, that is, the class of \poag s that are
isomorphic (as ordered groups) to some $\ZZ^n$, where $n\in\ZZ^+$.
Simplicial groups are basic building blocks of dimension groups, as
shows the following theorem, proved first, in a slightly different
form, in the context of commutative semigroups by Grillet, see
\cite{Gril76}, and then, in the context of \poag s, by Effros,
Handelman, and Shen, see \cite{EHS80}, and also \cite{Gpoag}:

\begin{theorem}\label{T:GEHS}
A \poag\ is a dimension group if and only if it is a direct limit of
simplicial groups and positive homomorphisms.\qed
\end{theorem}

In formula, if $G$ is a dimension group, then there exist a directed
set $I$, a direct system $\vv<S_i\mid i\in I>$ of simplicial groups,
together with transition homomorphisms $f_{ij}\colon S_i\to S_j$, for
$i\leq j$ in $I$, such that
 \begin{equation}\label{Eq:DirLimFla}
 G=\varinjlim_iS_i.
 \end{equation}
The following very simple example
shows that one cannot assume, in general, that the
$f_{ij}$ are \emph{positive embeddings}.

\begin{example}\label{E:Zalpha}
Let $G$ be an additive subgroup of $\RR$, containing as elements $1$
and an irrational number $\alpha$. We can assume that $\alpha>0$.

We claim that we cannot represent $G$ as a direct limit of simplicial
groups as in (\ref{Eq:DirLimFla}), with transition homomorphisms being
\emph{positive embeddings}. Indeed, suppose, to the contrary, that there
exists such a representation of
$G$. Let $i\in I$ such that $1$, $\alpha\in S_i$. Identifying $S_i$
with $\ZZ^n$ for some $n\in\ZZ^+$, write
 \[
 1=\vv<u_1,\ldots,u_n>\text{ and }\alpha=\vv<v_1,\ldots,v_n>.
 \]
Note that all the $u_i$ and $v_i$ belong to $\ZZ^+$. Define a rational
number $r$, by the formula
 \[
 r=\min\left\{\frac{v_i}{u_i}\mid i\in\fs(n)
 \text{ and }u_i>0\right\}.
 \]
Since the transition homomorphisms are positive embeddings, so are the
limiting maps $S_j\to G$, for $j\in I$. By taking $j=i$, we obtain that
 \[
 \{\vv<p,q>\in\ZZ\times\NN\mid p\leq q\alpha\}=
 \{\vv<p,q>\in\ZZ\times\NN\mid p\leq qr\}.
 \]
However, since $\alpha$ is irrational, this is impossible.
\end{example}

However, all the examples mentioned above can be expressed as direct
limits of simplicial groups and \emph{one-to-one} positive
homomorphisms. This is a special case of a much more general result,
due to Elliott, see \cite{Elli1,Elli2}, where it is proved that every
totally ordered abelian group can be expressed that way.
A more elaborate example, also due to Elliott, where the transition
maps $f_{ij}$ cannot all be taken to be one-to-one positive
homomorphisms, is presented in Remark 2.7 in \cite{Elli2}.

We see already that several kinds of dimension groups can be defined,
according to what kind of direct limit they are.

\begin{definition}\label{D:UltraS}
Let \eC\ be a class
of positive homomorphisms of \poag s (resp. \povs s). A dimension group
(resp. dimension vector space) is \emph{\eC-ultrasimplicial}, if it is a
direct limit of simplicial groups (resp. vector spaces) and
homomorphisms in \eC.
\end{definition}

Three classes of homomorphisms will be of particular importance to us.

\begin{itemize}
\item The class \eP\ of positive homomorphisms. By
Theorem~\ref{T:GEHS}, the class of
\eP-ultrasimplicial groups and
the class of dimension groups coincide.

\item The class \eI\ of \emph{one-to-one} positive homomorphisms.
The corresponding groups are thus called \eI-ultrasimplicial.
In \cite{Elli2}, they are called \emph{ultrasimplicial}.

\item The class \eE\ of \emph{positive embeddings}. The class of
\eE-ultrasimplicial groups will be of particular importance in this
paper.
\end{itemize}

Thus, Example~\ref{E:Zalpha} gives a \eI-ultrasimplicial group that
is not \eE-ultrasimplicial, while Elliott's example in \cite{Elli2}
is a dimension group that is not \eI-ultrasimplicial. This
already shows that none of the following implications
 \[
 \text{\eE-ultrasimplicial}\Rightarrow\text{\eI-ultrasimplicial}
 \Rightarrow\text{\eP-ultrasimplicial}
 \]
can be reversed.

Similarly, one can define \eP-\ (resp. \eI, \eE)-ultrasimplicial
\povs s over a given totally ordered division ring $K$. This time, \eP\
denotes the class of all positive \emph{linear maps} between vector
spaces over $K$, while \eI\
denotes the class of all one-to-one maps in \eP, and \eE\ denotes the
class of all order-embeddings in \eP.

Again, for a \povs\ $E$, $E$ is \eP-ultrasimplicial
if and only if it is a dimension vector space, see Theorem~7.2 in
\cite{Wehr}. It is an open problem whether, for \povs s,
\eP-ultrasimpliciality implies
\eI-ultrasimpliciality, even in the special case where $K=\QQ$.

The proof of the following lemma is an easy exercise.

\begin{lemma}\label{L:EI-ultra}
Let $G$ be a dimension group.
\begin{enumerate}
\item $G$ is \eI-ultrasimplicial if and only if every finite subset
of $G^+$ is contained into the submonoid generated by some finite
\emph{linearly independent} (over $\ZZ$) subset of $G^+$.

\item $G$ is \eE-ultrasimplicial if and only if every finite subset
of $G^+$ (resp. $G$) is contained into a simplicial subgroup of $G$.
\qed
\end{enumerate}
\end{lemma}

We leave to the reader the corresponding statements for \povs s.

Now let us recall the definition of a finitely presented \poag, as
given in \cite{Wehr}. As observed in this reference, it is equivalent
to the general definition of a finitely presented structure, given in
the language of universal algebra in \cite{Malc}.

\begin{definition}\label{D:FinPresPoag}
Let $G$ be a \poag. We say that $G$ is \emph{finitely presented}, if
there exists a finite generating subset $\{g_1,\ldots,g_n\}$ of $G$
such that the following set
 \[
 \{\vv<x_1,\ldots,x_n>\in\ZZ^n\mid
 g_1x_1+\cdots+g_nx_n\geq0\}
 \]
is a finitely generated submonoid of $\ZZ^n$.
\end{definition}

In particular, the validity of the definition above does not depend of
the choice of the finite generating subset $\{g_1,\ldots,g_n\}$, see
V.11, Corollary~7, p. 223 in \cite{Malc}; \emph{cf.} also \cite{Wehr}.

Observe that every simplicial group is finitely presented,
as a \poag. Hence, any \eE-ultrasimplicial group $G$ has the property
that every finite subset of $G$ is contained into a finitely presented
ordered subgroup of $G$. We shall state a slightly more convenient
equivalent form of this property, called \emph{coherence}, in
Corollary~\ref{C:EquivCohGrp}.

\section{Characterizations of finitely presented partially ordered
abelian groups}

We begin this section by stating a first characterization result of
finitely presented partially ordered abelian groups. It is a particular
case of a much more general result, see Theorem~8.1 in \cite{Wehr}.

\begin{proposition}\label{P:CharFPgroup}
Let $G$ be a \poag. Then $G$ is finitely presented if and only if $G$
is finitely generated as a group, and $G^+$ is finitely generated as
a monoid.\qed
\end{proposition}

We shall present here an alternative characterization of finitely
presented \poag s, based on purely order-theoretical properties of
$G$. To prepare for this characterization, let us first introduce some
terminology, and a lemma.

Let $\vv<M,+,0>$ be a \cm. We endow $M$ with the
\emph{algebraic} preordering $\leq$. By definition, for $x$, $y\in M$,
we define $x\leq y$ to hold, if there exists $z$ such that $x+z=y$. For
every subset $X$ of $M$, we put
 \[
 \upw X=\{y\in M\mid (\exists x\in X)(x\leq y)\}.
 \]
The origin of the argument of the following lemma can be traced back to
P. Freyd's very simple proof of Redei's Theorem, see \cite{Frey}, the
latter stating that every finitely generated \cm\ is finitely presented.
See also Proposition~12.7 in G. Brookfield's dissertation \cite{Broo},
where it is proved that every finitely generated cancellative
commutative monoid is well-founded,

\begin{lemma}\label{L:FinBasisMon}
Let $M$ be a finitely generated \cm. For every subset $X$ of $M$,
there exists a finite subset $Y$ of $X$ such that $\upw X=\upw Y$.
\end{lemma}

\begin{proof}
Let $\{g_1,\ldots,g_n\}$ be a finite generating subset of $M$.
Consider the polynomial ring $\ZZ[\xx_1,\ldots,\xx_n]$ on $n$
indeterminates on $\ZZ$.

We shall make use of the \emph{semigroup algebra} $\ZZ[M]$.
The elements of $\ZZ[M]$ are finite linear combinations
$\sum_{u\in M}n_u\dot u$,
where $\vv<n_u\mid u\in M>$ is an almost
null family of integers. The multiplication is defined in $\ZZ[M]$ in
such a way that if $w=u+v$ in $M$, then $\dot w=\dot u\cdot\dot v$.

There exists a unique ring homomorphism
$\varphi\colon\ZZ[\xx_1,\ldots,\xx_n]\to\ZZ[M]$ such
that $\varphi(\xx_i)=\dot g_i$, for all $i\in\fs(n)$. Note that
$\varphi$ is \emph{surjective}. By the Hilbert basis Theorem,
$\ZZ[\xx_1,\ldots,\xx_n]$ is n\oe therian; therefore, $\ZZ[M]$ is also
n\oe therian.

For every subset $X$ of $M$, let $I_X$ be the ideal of $\ZZ[M]$
generated by the subset
$\dot I=\{\dot u\mid u\in X\}$. It is clear that $I_X$
is the directed union of all the $I_Y$, where $Y$ is a finite subset of
$X$. Since $\ZZ[M]$ is n\oe therian, there exists a finite subset $Y$ of
$X$ such that $I_X=I_Y$. We prove that $\upw X=\upw Y$. It is clear
that $\upw X\supseteq\upw Y$. For the converse, it suffices to prove
that $X\subseteq\upw Y$. Let $x\in X$. Then $\dot x\in I_X=I_Y$, thus
there are elements $P_y$, $y\in Y$, of $\ZZ[M]$ such that
 \[
 \dot x=\sum_{y\in Y}P_y\dot y.
 \]
Thus, there exists $y\in Y$ such that the $\dot x$ component of
$P_y\dot y$ is nonzero. It follows easily that $y\leq x$; whence
$x\in\upw Y$.
\end{proof}

Recall that a partially ordered set $\vv<P,\leq>$ is
\emph{well-founded}, if every nonempty subset of $P$ admits a minimal
element. Equivalently, there exists no infinite strictly decreasing sequence
of elements of $P$.

If $\vv<P,\leq>$ is a partially ordered set and $X$ is a subset of
$P$, we denote by $\Min X$ the set of all minimal elements of
$X$. The following result gives equivalent condition, for a given \poag, to
be finitely presented.

\begin{theorem}\label{T:CharFPgroup}
Let $G$ be a \poag. Then the following are equivalent:
\begin{enumerate}
\item $G$ is finitely presented.

\item $G$ is finitely generated as a group, and $G^+$ is finitely
generated as a monoid.

\item $G$ is finitely generated as a group, $G^+$ is well-founded,
and, for every nonempty subset $X$ of $G^+$,
$\Min X$ is finite.

\item $G$ is finitely generated as a group, $G^+$ is well-founded, and
$\Min G^{++}$ is finite.
\end{enumerate}
\end{theorem}

\begin{proof}
The equivalence of (i) and (ii) has already been stated in
Proposition~\ref{P:CharFPgroup}.

(ii)$\Rightarrow$(iii) Assume that (ii) holds. Let $X$ be a nonempty
subset of $G^+$. By Lemma~\ref{L:FinBasisMon}, applied to the monoid
$G^+$, there exists a finite subset $Y$ of $X$ such that
$\upw X=\upw Y$. It follows that $\Min X=\Min Y$ is finite and nonempty.

(iii)$\Rightarrow$(iv) is trivial.

(iv)$\Rightarrow$(ii) Assume that (iv) holds. To conclude, it remains
to prove that $G^+$ is a finitely generated monoid.

\begin{sclaim}
Every nonempty majorized subset of $G$ admits a maximal element.
\end{sclaim}

\begin{scproof}
Let $X$ be a nonempty subset of $G$, majorized by $a\in G$. Since
$G^+$ is well-founded, the set $\{a-x\mid x\in X\}$ admits a minimal
element. This element has the form $a-x$, where $x$ is a maximal
element of $X$.
\end{scproof}

Now put $\Delta=\Min G^{++}$. By assumption, $\Delta$ is finite. Let
$M$ be the submonoid of $G^+$ generated by $\Delta$; we prove that
$M=G^+$. To prove the non trivial containment, let $a\in G^+$. Put
$X=\{x\in M\mid x\leq a\}$. Since $0\in X$, it follows from the Claim
above that $X$ admits a maximal element, say, $u$. Suppose that $u<a$.
Since $G^+$ is well-founded, there exists $v\in\Delta$ such that
$v\leq a-u$. Hence $u+v\in X$, which contradicts the maximality of $u$
in $X$. Hence $u=a$, which proves that $a\in M$. Thus we have proved
that $G^+$ is generated by the finite subset $\Delta$.
\end{proof}

The following examples, \ref{E:nonwf} to \ref{E:nonfinmin}, show the
independence of the conditions listed in (iv) of
Theorem~\ref{T:CharFPgroup}.

\begin{example}\label{E:nonwf}
This example is a particular case of Example~\ref{E:Zalpha}. Let
$\alpha$ be a real irrational number, and put $G=\ZZ+\ZZ\alpha$,
viewed as an ordered additive subgroup of $\RR$. Then $G$ is a
finitely generated \poag, $\Min G^{++}=\varnothing$ is finite, but
$G^+$ is not well-founded.
\end{example}

\begin{example}
Let $G$ be any non finitely generated abelian group (for example, the
free abelian group on any infinite set), endowed with the trivial
positive cone $G^+=\{0\}$. Then $G^+$ is well-founded and
$\Min G^{++}$ is finite, but $G$ is not finitely generated.
\end{example}

\begin{example}\label{E:nonfinmin}
Endow $G=\ZZ\times\ZZ$ with the positive cone defined by
$G^+=\{\vv<0,0>\}\cup(\NN\times\NN)$. Then $G$ is finitely generated,
$G^+$ is well-founded, but $\Min G^{++}$ consists of all elements of
the form $\{\vv<1,n>\}$ or $\{\vv<n,1>\}$, where $n\in\NN$, thus it is
infinite.
\end{example}

The following result is obtained, in the much more general context
of partially ordered right modules over \emph{coherent} normally
ordered rings, and by completely different methods, in \cite{Wehr}.

\begin{corollary}\label{C:Subgrp}
Let $G$ be a finitely presented \poag. Then every ordered subgroup of
$G$ is finitely presented.
\end{corollary}

\begin{proof}
By the elementary theory of abelian groups, every subgroup of a
finitely generated abelian group is finitely generated. We conclude
by using the characterization (iii) of Theorem~\ref{T:CharFPgroup} for
finitely presented \poag s.
\end{proof}

Note the following useful corollary of Theorem~\ref{T:CharFPgroup}:

\begin{corollary}\label{C:Subgrp2}
Let $G$ be a finitely presented \poag. Then $G$ is \emph{weakly
Archimedean}, that is, the condition
 \begin{equation}\label{Eq:warch}
 (\forall n\in\NN)(na\leq b)
 \end{equation}
implies that $a=0$, for all $a$, $b\in G^+$.
\end{corollary}

\begin{proof}
Suppose that (\ref{Eq:warch}) is satisfied, with $a>0$. Let $H$ be the
lexicographical product of $\ZZ$ with itself; thus, we have
 \[
 H=\ZZ\times\ZZ,\quad\text{and}\quad
 H^+=\{\vv<m,n>\in\ZZ^+\times\ZZ\mid m=0\Rightarrow n\geq0\}.
 \]
Then the map $\vv<m,n>\mapsto mb+na$ is an order-embedding from $H$
into $G$. Since $G$ is finitely presented, it follows from
Corollary~\ref{C:Subgrp} that $H$ is finitely presented. By
Theorem~\ref{T:CharFPgroup}, it follows that $H^+$ is a finitely
generated monoid, which is easily verified not to be the case.
\end{proof}

\begin{corollary}\label{C:FinInt}
Let $G$ be a finitely presented \poag. For all elements $a$ and $b$ of
$G$ such that $a\leq b$, the \emph{interval} $[a,b]$, defined by
 \[
 [a,b]=\{x\in G\mid a\leq x\leq b\},
 \]
is finite.
\end{corollary}

\begin{proof}
Without loss of generality, $a=0$. By Theorem~\ref{T:CharFPgroup},
$G^+$ admits a finite generating subset, say, $\{g_1,\ldots,g_m\}$,
where $g_i>0$ for all $i$. By Corollary~\ref{C:Subgrp2}, for all
$i\in\fs(m)$, there exists a largest integer $n_i\in\ZZ^+$ such that
$n_ig_i\leq b$. Put $n=\max\{n_1,\ldots,n_m\}$.
Let $x\in G^+$ such that $0\leq x\leq b$. Since
$x\in G^+$, there are non-negative integers $k_1$, \dots, $k_m$ such
that $x=\sum_{i=1}^mk_ig_i$.
The inequality $k_ig_i\leq b$ holds for all $i\in\fs(m)$, thus
$k_i\leq n_i$. Since
$k_i\geq0$, one can thus define a surjection from
$\prod_{i=1}^m\{0,1,\ldots,n_i\}$ to
$[0,b]$, by the rule $\vv<k_1,\ldots,k_m>\mapsto\sum_{i=1}^mk_ig_i$; in
particular, $[0,b]$ is finite.
\end{proof}

\begin{remark}
By using Theorem~\ref{T:DescFPvs}, one can prove a similar result for
a \emph{\povs} $E$ over $\RR$. Of course, the conclusion does not
state any longer that intervals of the form $[a,b]$ are finite, but,
rather, that they are \emph{compact} for the natural topology of $E$.
\end{remark}

On the other hand, the results of Corollaries \ref{C:FinInt} and
\ref{C:Subgrp2} do not give a characterization of finitely presented
\poag s, as shows Example~\ref{E:nonfinmin}.

\begin{example}\label{E:NotArch}
Define a \poag\ to be \emph{Archimedean}, if, for all $a$, $b\in G$,
if $na\leq b$ for all $n\in\ZZ^+$, then $a\leq0$. Thus every
Archimedean \poag\ is weakly Archimedean. It is also a classical
result that a directed partially ordered group is Archimedean if and
only if it can be embedded into a Dedekind complete lattice-ordered
group, see \cite{BKWo}. In particular, every Archimedean partially
ordered group is \emph{commutative} and \emph{unperforated}.

We present here an example of directed, torsion-free, non
Archimedean, finitely presented \poag\ $G$. Define it by
 \[
 G=\ZZ\times\ZZ,\quad\text{and}\quad
 G^+=\{\vv<0,0>\}\cup\{\vv<x,y>\in\ZZ^+\times\ZZ^+\mid x+y\geq2\}.
 \]
Then $G$ is finitely generated as a group, and $G^+$ is finitely
generated as a monoid (by the seven elements $\vv<2,0>$, $\vv<1,1>$,
$\vv<0,2>$, $\vv<3,0>$, $\vv<2,1>$, $\vv<1,2>$, $\vv<0,3>$). Thus $G$
is a finitely presented \poag.

However, $G$ is not unperforated. For example, if we take
$a=\vv<1,0>$, then $G$ satisfies $2a\geq0$, but $a\ngeq0$. It follows
that $G$ is not Archimedean. This can also be verified directly.
Indeed, put $b=\vv<2,0>$. Then the inequality
$n(-a)\leq b$ holds for all $n\in\ZZ^+$, but $a\ngeq0$.
\end{example}

The following definition is a particular case of a general notion
studied in \cite{Wehr}.

\begin{definition}\label{D:CohOrdGrp}
Let $G$ be a \poag\ (resp. a \povs). We say that $G$ is \emph{coherent},
if every finitely generated ordered subgroup (resp. subspace) of $G$ is
finitely presented.
\end{definition}

For a given \poag, several equivalent definitions for being finitely
presented are stated in \cite{Wehr}. Among those, note the following:
\emph{A \poag\ $G$ is coherent, if and only if
the solution set of every homogeneous system of $m$ equations and
inequalities with parameters from $G$ and $n$ unknowns in $\ZZ$ is a
finitely generated submonoid of $\ZZ^n$, for all $m$, $n\in\NN$}. In view
of Corollary~\ref{C:Subgrp}, we can state immediately the following
corollary:

\begin{corollary}\label{C:EquivCohGrp}
Let $G$ be a \poag. Then $G$ is coherent if and only if every finite
subset of $G$ is contained into a finitely presented ordered subgroup
of $G$.\qed
\end{corollary}

Note the following corollary about \eE-ultrasimplicial groups:

\begin{corollary}\label{C:UltrFP}
Every \eE-ultrasimplicial dimension group is coherent.\qed
\end{corollary}

\begin{example}
Let $I$ be any set. Then the free abelian group $F_I=\ZZ^{(I)}$
on $I$, endowed with the positive cone $(\ZZ^+)^{(I)}$, is a coherent
\poag\ (because it is \eE-ultrasimplicial). However, if $I$ is
infinite, then $F_I$ is not finitely presented, because it is not
finitely generated as a group.
\end{example}

\section{Convex domains over totally ordered division rings}

We shall fix in this section a totally ordered division ring $K$.
We shall restate, in the context of right vector spaces over $K$, some
elementary theory of convex polyhedra, with, for example, a
corresponding version of the Hahn-Banach Theorem. Of course, in the
context of real vector spaces, these results are plain classical
mathematics. However, we have been unable to locate a reference where
they are proved in the present more general context (to which the usual
topological proofs do not extend), thus we provide the proofs here.

If $E$ is a right vector space over $K$, an \emph{affine functional}
over $E$ is a map of the form $x\mapsto p(x)+b$, where
$p\colon E\to K$ is a linear functional on $E$ and $b$ is a
fixed element of $K$. If $E=K^n$ is given its natural structure of
right vector space over $K$, then the affine functionals on $E$ are
exactly the maps from $E$ to $K$ of the form
 \[
 \vv<\xi_1,\ldots,\xi_n>\mapsto a_1\xi_1+\cdots+a_n\xi_n+b,
 \]
where $a_1$, \dots, $a_n$, $b$ are fixed elements of $K$.
A subset of $K^n$ is \emph{bounded}, if it is contained into a product
of bounded intervals of $K$. Observe that this notion does not depend
on the chosen basis of the ground vector space.

If $p$ is an affine functional on $E$, then we put
 \[
 \|p\geq0\|=\{x\in E\mid p(x)\geq0\},
 \]
the \emph{upper half space} defined by $p$. A \emph{convex domain} is
a finite intersection of half spaces. A \emph{convex polyhedron} is
any \emph{bounded} convex domain. It is straightforward to verify that
any convex domain $U$ is, indeed, \emph{convex}, that is,
$x(1-\lambda)+y\lambda\in U$, for all $x$,
$y\in U$ and all $\lambda\in K$ such that $0\leq\lambda\leq1$,.

The following result is an immediate consequence of the definition of
a convex domain.

\begin{proposition}[Hahn-Banach Theorem]
\label{P:HB}
Let $E$ be a right vector space over $K$, let $C$ be a convex domain
of $E$ such that $0\notin C$. Then there exists a linear functional
$p$ on $E$ such that $p|_C\geq1$.
\end{proposition}

\begin{proof}
Since $C$ is the intersection of half spaces, there exists an affine
functional $q$ on $E$ such that $C\subseteq \| q\geq0 \|$ but $q(0)<0$.
Put $p=(-1/q(0))(q-q(0))$.
\end{proof}

We record a particular case of this statement, that will prove useful
in the proof of Theorem~\ref{T:EUltr}.

\begin{lemma}\label{L:HyperplHB}
Let $E$ be a right vector space over $K$, let $C$ be a convex domain
of $E$, let $H$ be an affine hyperplane of $E$ such
that $C\subseteq H$ and $0\notin H$. Let $a$ be an element of
$H\setminus C$. Then there exists a linear functional $p$ on $E$ such
that $p(a)=0$ and $p|_C\geq1$.
\end{lemma}

\begin{proof}
Since $0\notin H$, there exists a linear functional $r$ on $E$ such
that $H=\{x\in E\mid r(x)=1\}$. Since $C$ is a convex domain and
since $a\notin C$, there exists an affine functional $q$ such that
$q(a)=0$ and $q|_C\geq1$. Define an affine functional $p$ on
$E$ by
 \[
 p(x)=q(x)+q(0)(r(x)-1),\qquad\text{for all }x\in E.
 \]
Since $p(0)=0$, $p$ is a linear functional on $E$. Furthermore,
$p|_H=q|_H$, which completes the proof.
\end{proof}

We shall now, in the finite dimensional case, characterize convex
polyhedra of $E$ as convex hulls of finite subsets of $E$.

\begin{lemma}\label{L:1ProjConv}
Let $E$ be a vector space over $K$, let $U$ be a convex domain of
$E\times K$. Then the projection of $U$ on $E$ is a convex domain of
$E$.
\end{lemma}

\begin{proof}
Let $p_1$, \dots, $p_n$ be affine functionals on $E\times K$ such that
 \[
 U=\bigcap_{i=1}^n\|p_i\geq0\|.
 \]
For all $i\in\fs(n)$, there exist an affine functional $q_i$ on $E$
and an element $\lambda_i$ of $K$ such that
$p_i(\vv<x,\xi>)=q_i(x)+\lambda_i\xi$, for all
$\vv<x,\xi>\in E\times K$.

Put $X=\{i\mid\lambda_i=0\}$, $Y=\{i\mid\lambda_i>0\}$ and
$Z=\{i\mid\lambda_i<0\}$.

Denote by $V$ the projection of $U$ on $E$. By definition, $V$ is the
set of all elements $x$ of $E$ such that there exists $\xi\in K$
satisfying the following conditions:
 \begin{align*}
 0&\leq q_i(x),&\text{for all }i\in X,\\
 \xi&\geq -\frac{1}{\lambda_j}q_j(x), & \text{for all }j\in Y,\\
 \xi&\leq -\frac{1}{\lambda_k}q_k(x), & \text{for all }k\in Z.
 \end{align*}
Therefore, an element $x$ of $E$ belongs to $V$ if and only if $x$
satisfies the following conditions:
 \begin{align*}
 0&\leq q_i(x),&\text{for all }i\in X,\\
 -\frac{1}{\lambda_j}q_j(x)&\leq-\frac{1}{\lambda_k}q_k(x),&
 \text{for all }\vv<j,k>\in Y\times Z.
 \end{align*}
Hence, $V$ is the intersection of the half spaces $\|p_i\geq0\|$, for
$i\in X$, and\linebreak
$\|(1/\lambda_j)q_j-(1/\lambda_k)q_k\geq0\|$, for
$\vv<j,k>\in Y\times Z$. Hence, $V$ is a convex domain of $E$.
\end{proof}

Thus, an easy induction on the dimension of $F$ yields the following
corollary:

\begin{corollary}\label{C:ProjConv}
Let $E$ and $F$ be right vector spaces over $K$, with $F$ finite
dimensional, and let $U$ be a convex domain of $E\times F$. Then the
projection of $U$ on $E$ is a convex domain of $E$.\qed
\end{corollary}

For any subset $X$ of a right vector space $E$ over $K$, define
$\Conv X$ to be the set of all \emph{convex combinations} of elements
of $X$, that is, the set of all elements of $E$ of the form
$\sum_{i=1}^nx_i\lambda_i$, where $n\in\NN$, the $x_i$ belong to $X$,
the $\lambda_i$ belong to $K^+$ and $\sum_{i=1}^n\lambda_i=1$.

\begin{proposition}\label{P:ConvFin}
Let $E$ be a finite dimensional right vector space over $K$, let $X$ be
a finite subset of $E$. Then $\Conv X$ is a convex polyhedron of $E$.
\end{proposition}

\begin{proof}
Without loss of generality, $E=K^m$, for some $m\in\ZZ^+$.

It is obvious that $\Conv X$ is bounded.

Write $X=\{a_1,\ldots,a_n\}$, and let $Y$ be the subset of $E\times K^n$
defined by
 \[
 Y=\left\{\vv<x,\vv<\xi_1,\ldots,\xi_n>>\in E\times(K^+)^n\mid
 \sum_{j=1}^n\xi_j=1\text{ and }x=\sum_{j=1}^na_j\xi_j\right\}.
 \]
In particular, $\Conv X$ is the projection of $Y$ on $E$.
We write $a_j=\vv<a_{1j},\ldots,a_{mj}>$, for all $j\in\fs(n)$.
An element $\vv<x,\vv<\xi_1,\ldots,\xi_n>>$ of $E\times K^n$, with
$x=\vv<x_1,\ldots,x_m>$, belongs to $Y$ if and only if the following
conditions hold:
 \begin{equation*}
 \begin{cases}
 x_i-\sum_{j=1}^na_{ij}\xi_j=0,&\text{for all }i\in\fs(m),\\
 \hfill\xi_j\geq0,&\text{for all }j\in\fs(n),\\
 \hfill 1-\xi_j\geq0,&\text{for all }j\in\fs(n),\\
 \hfill\sum_{j=1}^n\xi_j-1=0.&
 \end{cases}
 \end{equation*}
Since $\xi=0$ if and only if $\xi\geq0$ and $-\xi\geq0$, we obtain
that $Y$ is a convex domain of $E\times K^m$. By
Corollary~\ref{C:ProjConv}, the projection of $Y$ on $E$ is a convex
domain of $E$.
\end{proof}

In particular, we recover the fact that $\Conv X$ is the
\emph{convex hull} of $X$, that is, the least convex subset of $E$
containing $X$. Of course, this fact is, also, much easier to
establish directly.

\begin{corollary}\label{C:CharPolyh}
Let $E$ be a finite dimensional right vector space over $K$. Then the
convex polyhedra of $E$ are exactly the subsets of the form $\Conv X$,
where $X$ is a finite subset of $E$.
\end{corollary}

\begin{proof}
One direction follows from Proposition~\ref{P:ConvFin}. To prove the
other direction, we prove, by induction on the dimension of $E$, that
every convex polyhedron of $E$ has the form $\Conv X$, for some finite
subset $X$ of $E$.
 It is trivial if $\dim E=0$. Suppose, now, that
$\dim E>0$, and let $C$ be a convex polyhedron of $E$. By definition,
there exist $n\in\NN$ and nonzero affine functionals $p_j$,
$j\in\fs(n)$, such that
 \[
 C=\bigcap_{j=1}^n\|p_j\geq0\|.
 \]
For all $i\in\fs(n)$, define a convex domain $C_i$ by
 \[
 C_i=\|p_i=0\|\cap\bigcap_{j\ne i}\|p_j\geq0\|,
 \]
where we write $\|p_i=0\|=\{x\mid p_i(x)=0\}$. In particular,
$C_i$ is contained into $C$, thus $C_i$ is a convex polyhedron.
Furthermore, $C_i$ is contained into the affine hyperplane
$H_i=\|p_i=0\|$, thus, by induction hypothesis, there exists a finite
subset $X_i$ of
$H_i$ such that $C_i=\Conv X_i$. Thus, to conclude, it suffices to prove
that $C$ is contained into the convex hull of $\bigcup_{i=1}^nC_i$,
because it follows then immediately that $C$ is equal to the convex
hull of the finite set $\bigcup_{i=1}^nX_i$.
Thus, let $x\in C$. If $p_i(x)=0$ for some $i$, then we are done,
because $x$ belongs to $C_i$; so, suppose that $p_i(x)>0$ for all $i$.
Let $u$ be an arbitrary element of $E\setminus\{0\}$. We put
$q_i=p_i-p_i(0)$, for all $i\in\fs(n)$. We define subsets $I_+$ and $I_-$
of
$\fs(n)$ as follows:
 \begin{align*}
 I_+&=\{i\in\fs(n)\mid q_i(u)>0\},\\
 I_-&=\{i\in\fs(n)\mid q_i(u)<0\}.
 \end{align*}
If $I_+=\varnothing$, then
$p_i(x-u\lambda)=p_i(x)-q_i(u)\lambda\geq p_i(x)>0$, for all
$\lambda\in K^+$ and all $i\in\fs(n)$, thus
$x-u\lambda\in C$; this contradicts the fact that $C$ is bounded.
Thus, $I_+$ is nonempty. Similarly, $I_-$ is nonempty. Define two
elements $\alpha$ and $\beta$ of $K^{++}$ as follows:
 \begin{align*}
 \alpha&=\min\left\{
 q_i(u)^{-1}p_i(x)\mid i\in I_+\right\},\\
 \beta&=\min\left\{
 -q_i(u)^{-1}p_i(x)\mid i\in I_-\right\}.
 \end{align*}
Let $j\in I_+$ and $k\in I_-$ such that $\alpha=q_j(u)^{-1}p_j(x)$ and
$\beta=-q_k(u)^{-1}p_k(x)$. Then
$x-u\alpha$ belongs to $C_j$ and $x+u\beta$ belongs to $C_k$.
Moreover, we have
 \[
 x(1+\beta^{-1}\alpha)=(x-u\alpha)1+(x+u\beta)\beta^{-1}\alpha,
 \]
so that $x$ belongs to the convex hull of $C_j\cup C_k$.
\end{proof}

\section{Finitely presented partially ordered vector spaces}
\label{S:FPvs}

In order to give a representation of finitely presented partially
ordered vector spaces (Theorem~\ref{T:DescFPvs}), we first prove a
lemma.

\begin{lemma}\label{L:Cpos}
Let $E$ be a finite dimensional right vector space over $K$, let $C$
be a convex polyhedron of $E$ such that $0\notin C$. Put
$n=\dim E$, and endow $K^n$ with its natural structure of
\povs\ over $K$.

Then there exists an isomorphism of vector spaces,
$\varphi\colon E\to K^n$, such that
$\varphi[C]\subseteq(K^{++})^n$.
\end{lemma}

\begin{proof}
If $C=\varnothing$, it is trivial. Thus suppose that $C$ is nonempty.
In particular, $n>0$. Thus put $n=m+1$, where $m\in\ZZ^+$.

By Proposition~\ref{P:HB}, there exists a linear functional $p$ on $E$
such that $p|_C\geq1$. Put $H=\ker p$, and let
$\vv<e_i\mid 0\leq i<m>$ be a basis of $H$. Furthermore, let $e\in E$
such that $p(e)=1$. Let $X$ be a finite subset of $E$ such that
$C=\Conv X$. For all $x\in X$, we decompose $x$ is the basis formed by
the $e_i$ and $e$, as follows:
 \begin{equation}\label{Eq:Exprx}
 x=\sum_{0\leq i<m}e_i\cdot\xi_{x,i}+e\cdot p(x).
 \end{equation}
Note, in particular, that $p(x)\geq 1$.
Let $\lambda$ be an element of $K$ such that
$\lambda>-\xi_{x,i}p(x)^{-1}$,
for all $x\in X$ and all $i<m$.
Put $f=e-\left(\sum_{0\leq i<m}e_i\right)\cdot\lambda$. Let $\varphi$
the unique linear map from $E$ to $K^n$ sending the basis
$\vv<e_0,\ldots,e_{m-1},f>$ of $E$ to the canonical basis of $K^n$.
Since $e=f+\left(\sum_{0\leq i<m}e_i\right)\cdot\lambda$, it follows
from (\ref{Eq:Exprx}) that
\begin{align*}
 x&=\sum_{0\leq i<m}e_i\cdot\xi_{x,i}+
 \left[f+\left(\sum_{0\leq i<m}e_i\right)\cdot\lambda\right]p(x)\\
 &=\sum_{0\leq i<m}e_i\cdot(\xi_{x,i}+\lambda p(x))+f\cdot p(x),
 \end{align*}
for all $x\in X$,
thus all components of $x$ in the basis $\vv<e_0,\ldots,e_{m-1},f>$
belong to $K^{++}$. This means that $\varphi(x)\in(K^{++})^n$, for
all $x\in X$. Hence the same conclusion holds for all $x\in C$.
\end{proof}

Recall that if $E$ is a right vector space over a totally ordered
division ring $K$, we say that a \emph{$K^+$-subsemimodule} of $E$ is a
submonoid $M$ of $E$ such that $K^+M=M$.

We also record the following fact, which has been established in
\cite{Wehr}:

\begin{proposition}\label{P:Kcoh}
Let $K$ be a totally ordered division ring, and let $E$ be a \povs\ over
$K$. Then the following properties hold:
\begin{enumerate}
\item $E$ is finitely presented if and only if $E$ is a finite
dimensional $K$-vector space and $E^+$ is a finitely generated
$K^+$-semimodule.

\item Suppose that $E$ is finitely presented. Then every subspace of
$E$, endowed with the induced ordering, is finitely presented.\qed
\end{enumerate}
\end{proposition}

Observe that even for $K=\QQ$, the analogue of
Theorem~\ref{T:CharFPgroup} does not hold for \povs s over
$K$; indeed, $\QQ$, endowed with its natural ordering, is
a finitely presented \povs\ over $\QQ$, but
$\QQ^+$ is not well-founded. Hence, Theorem~\ref{T:CharFPgroup} is
something extremely specific to ordered groups. In particular, one
cannot establish Proposition~\ref{P:Kcoh}(ii) by using an analogue of
the proof of Corollary~\ref{C:Subgrp}.

\begin{theorem}\label{T:DescFPvs}
Let $K$ be a totally ordered division ring. Then the finitely
presented \povs s over $K$ are, up to isomorphism, those of the form
$\vv<K^n,+,0,P>$, where $n\in\ZZ^+$ and $P$ is a finitely generated
$K^+$-subsemimodule of $(K^+)^n$.
\end{theorem}

\begin{remark}
As the proof will show, one can refine $(K^+)^n$ into
$\{0\}\cup(K^{++})^n$.
\end{remark}

\begin{proof}
Let $E$ be a \povs\ over $K$. If $E$ has
the indicated form, then, by Proposition~\ref{P:Kcoh}, $E$ is finitely
presented.

Conversely, suppose that $E$ is finitely presented. By
Proposition~\ref{P:Kcoh}, $E$ is finite dimensional, and $E^+$ is a
finitely generated $K^+$-subsemimodule of $E$. Put $n=\dim E$, and let
$X$ be a finite generating subset of $E^+$ such that $0\notin X$. Let
$C=\Conv X$ be the convex hull of $X$.
Since $E^+\cap(-E^+)=\{0\}$, $0$ does not belong to $C$. Thus, by
Lemma~\ref{L:Cpos}, there exists an isomorphism of vector spaces,
$\varphi\colon E\to K^n$, such that $\varphi[C]\subseteq(K^+)^n$. Take
$P=\varphi[E^+]$. Note that $P$ is a finitely generated
$K^+$-subsemimodule of $K^n$.
\end{proof}

\section{Finitely presented torsion-free partially ordered groups}

In this section, we shall obtain representation results for
finitely presented \poag s, similar to those obtained in
Section~\ref{S:FPvs} for partially ordered vector spaces.

We shall need the following version of Lemma~\ref{L:Cpos}.

\begin{lemma}\label{L:Cposgroup}
Let $n$ be a natural number, let $X$ be a finite subset of $\ZZ^n$
such that $0$ is not a linear combination with coefficients in $\NN$
of a nonempty subset of $X$.

Then there exists a group automorphism $\varphi$ of $\ZZ^n$ such that
$\varphi[C]\subseteq(\ZZ^{++})^n$.
\end{lemma}

Note that the condition on $X$ can be expressed differently, by
saying that $0$ does not belong to the convex hull of $X$ in $\QQ^n$.

\begin{proof}
The proof is similar to the proof of Lemma~\ref{L:Cpos}, with a few
modifications in order to ensure that $\ZZ^n$ is invariant under the
automorphism $\varphi$ of $\QQ^n$ which is constructed. We just
indicate here those modifications. Since the result is, again, trivial
for $X=\varnothing$, we suppose that $X$ is nonempty. Thus, $n$ is
positive; put $m=n-1$.

First, of course, we embed $\ZZ^n$ into $\QQ^n$. As in the proof of
Lemma~\ref{L:Cpos}, there exist a linear functional $p$ on $\QQ^n$
such that $p|_X\geq1$. There are rational numbers $r_1$, \dots, $r_n$
such that
 \[
 p(\vv<x_1,\ldots,x_n>)=r_1x_1+\cdots+r_nx_n,
 \]
for all $x_1$, \dots, $x_n\in\QQ$.
After having multiplied $p$ by some suitable positive integer, one may
suppose that the $r_i$ are \emph{integers}. Furthermore, after
dividing the $r_i$ by their greatest common divisor, one may suppose
that $r_1$, \dots, $r_n$ are \emph{coprime}. By Bezout's Theorem,
there are integers $u_1$, \dots, $u_n$ such that
$\sum_{i=1}^nr_iu_i=1$. Put $e=\vv<u_1,\ldots,u_n>$.

Furthermore, $H=\ker p\cap\ZZ^n$ is a subgroup of rank $n-1$ of $\ZZ^n$.
In particular, it is free abelian. Therefore, $H$ admits a basis
$\vv<e_0,\ldots,e_{m-1}>$ over $\ZZ$.

The rest of the proof follows the pattern of the corresponding part of
the proof of Lemma~\ref{L:Cpos}, with the obvious modifications. For
example, $\lambda$ is, now, a large enough \emph{integer}.
One needs also to observe that both
$\vv<e_0,\ldots,e_{m-1},e>$ and $\vv<e_0,\ldots,e_{m-1},f>$ are,
indeed, generating subsets of $\ZZ^n$ (they are obviously independent
over $\ZZ$), thus the map $\varphi$ is a group automorphism of $\ZZ^n$.
\end{proof}

\begin{theorem}\label{T:DescFPgp}
The finitely presented, torsion-free \poag s are, up to
isomorphism, exactly those of the form $\vv<\ZZ^n,+,0,P>$, where
$P$ is a finitely generated submonoid of $(\ZZ^+)^n$.
\end{theorem}

\begin{proof}
Since $G$ is a finitely generated torsion-free abelian group, it is
isomorphic, as a group, to $\ZZ^n$, for some $n\in\ZZ^+$. Thus,
without loss of generality, $G=\ZZ^n$ as an abelian group. By
Theorem~\ref{T:CharFPgroup}, the additive monoid $G^+$ admits a finite
generating subset, $X$. One can of course suppose that $0\notin X$.
Then $0$ does not belong to the convex hull of
$X$ in $\QQ^n$, thus, by Lemma~\ref{L:Cposgroup}, there
exists a group automorphism $\varphi$ of $\ZZ^n$ such that
$\varphi[X]\subseteq(\ZZ^+)^n$.
Therefore, we also have $\varphi[G^+]\subseteq(\ZZ^+)^n$.

Hence, $G$ is, as a \poag, isomorphic to $\vv<\ZZ^n,+,0,\varphi[G^+]>$.
\end{proof}

\begin{remark}
It is well known that every submonoid of $\ZZ^+$ is finitely generated.
However, the analogue of this result does not hold for all higher
dimensions, as shows, for example, the submonoid of $\ZZ^+\times\ZZ^+$
generated by all pairs $\vv<n,1>$, where $n\in\ZZ^+$. This justifies
the additional precision in the statement of Theorem~\ref{T:DescFPgp},
that $P$ is finitely generated.
\end{remark}

We illustrate Theorem~\ref{T:DescFPgp} with two examples.

\begin{example}\label{E:a+b=2c}
Let $G$ (resp. $G^+$) be the subgroup (resp. submonoid) of $\ZZ^3$
generated by $a=\vv<2,0,1>$, $b=\vv<0,2,1>$, and $c=\vv<1,1,1>$. By
Theorem~\ref{T:CharFPgroup}, $G$ is a finitely presented \poag. In fact,
it is not hard to prove that
$G$ is the \poag\ defined by generators $a$, $b$, and $c$, and relations
 \[
 a\geq0;\quad b\geq0;\quad c\geq0;\quad a+b=2c.
 \]
Note that $G$ is a free abelian group of rank $2$
(with the two generators $a$ and $c$).
Thus, by Theorem~\ref{T:DescFPgp}, $G$ can be represented as $\ZZ^2$,
endowed with a finitely generated submonoid of $(\ZZ^+)^2$ as positive
cone. Here is such a representation. Take $a'=\vv<1,0>$, $b'=\vv<1,2>$
and $c'=\vv<1,1>$, and let $P$ be the submonoid of $(\ZZ^+)^2$
generated by $\{a',b',c'\}$. Then it is not hard to verify that
 \[
 \vv<G,+,0,G^+>\cong\vv<\ZZ^2,+,0,P>.
 \]
\end{example}

Our next example shows that there is no analogue of
Theorem~\ref{T:DescFPgp} for non finitely presented \poag s.

\begin{example}\label{E:Zalpha2}
Let $G$ be the \poag\ given in Example~\ref{E:nonwf}. Then $G$ is a
finitely generated \poag\ (it has two generators, $1$ and $\alpha$).
However, if $P$ is a submonoid of $(\ZZ^+)^2$, then $G$ cannot be
isomorphic to $\vv<\ZZ^2,+,0,P>$, because $G$ is totally ordered, while
one cannot have $P\cup(-P)=\ZZ^2$.
\end{example}

For our next application, we need the following folklore lemma.

\begin{lemma}\label{L:SysIneq}
Let $m$, $n\in\ZZ^+$, and let $a_{i,j}$ ($1\leq i\leq m$ and
$1\leq j\leq n$), and $b_i$ ($1\leq i\leq m$) be rational numbers.
Let $\Sigma$ be the following system of inequalities:
 \begin{equation}\label{Eq:GenSys}
 \begin{cases}
 a_{1,1}x_1+\cdots+a_{1,n}x_n&\leq b_1\\
 a_{2,1}x_1+\cdots+a_{2,n}x_n&\leq b_2\\
 \hfill\vdots\hfill&\hfill\vdots\hfill\\
 a_{m,1}x_1+\cdots+a_{m,n}x_n&\leq b_m
 \end{cases}
 \end{equation}

If $\Sigma$ admits a solution in $\RR^n$, then it admits a solution in
$\QQ^n$.
\end{lemma}

\begin{proof}
By induction on $n$. If $n=0$, then it is trivial. Suppose that
$n>0$ and that the result is proved for $n-1$. Note that
(\ref{Eq:GenSys}) can be rewritten under the following form:
 \begin{equation}\label{Eq:GenSys2}
 \begin{cases}
 a_{1,n}x_n&\leq c_1(\vec x)\\
 a_{2,n}x_n&\leq c_2(\vec x)\\
 \hfill\vdots\hfill&\hfill\vdots\hfill\\
 a_{m,n}x_n&\leq c_m(\vec x),
 \end{cases}
 \end{equation}
where we put $c_i(\vec x)=b_i-\sum_{1\leq j<n}a_{i,j}x_j$. Define
 \[
 I=\{i\mid a_{i,n}=0\},\quad U=\{i\mid a_{i,n}>0\},\quad
 \text{and}\quad V=\{i\mid a_{i,n}<0\}.
 \]
Then put $d_i(\vec x)=a_{i,n}^{-1}c_i(\vec x)$, for all $i\in U\cup V$.
Then the system (\ref{Eq:GenSys2}) can be rewritten as follows:
 \begin{equation}\label{Eq:GenSys3}
 \begin{cases}
 \hfil 0\leq c_i(\vec x)&\quad\text{for all }i\in I,\\
 x_n\leq d_u(\vec x)&\quad\text{for all }u\in U,\\
 x_n\geq d_v(\vec x)&\quad\text{for all }v\in V.
 \end{cases}
 \end{equation}
Whatever ground field $\QQ$ or $\RR$ we consider, the existence of a
solution of (\ref{Eq:GenSys3}) is equivalent to the existence of a
solution to the following system:
 \begin{equation*}
 \begin{cases}
 0\leq c_i(\vec x)&\quad\text{for all }i\in I,\\
 d_v(\vec x)\leq d_u(\vec x)&
 \quad\text{for all }\vv<u,v>\in U\times V,
 \end{cases}
 \end{equation*}
which is a system of the form (\ref{Eq:GenSys}), but with the $n-1$
unknowns $x_1$, \dots, $x_{n-1}$.
\end{proof}

It is to be noted that Lemma~\ref{L:SysIneq} remains valid for systems
containing \emph{strict} inequalities (with a similar proof), but
we shall not need this fact here.

\begin{corollary}\label{C:NonPerfAr}
Every finitely presented, unperforated \poag\ is Archimedean.
\end{corollary}

\begin{proof}
Let $G$ be a finitely presented, unperforated \poag. We prove that $G$
is Archimedean. Note, in particular, that $G$ is torsion-free. If
$G^+=\{0\}$, then the result is trivial, so suppose that $G^+\ne\{0\}$.
By Theorem~\ref{T:DescFPgp}, we can suppose that $G$ has the form
$\vv<\ZZ^m,+,0,G^+>$, where $G^+$ is a finitely generated submonoid of
$(\ZZ^+)^m$. Let $\{g_1,\ldots,g_n\}$ (with $n>0$) be a finite
generating subset of $G^+$, with all the $g_j$ nonzero. We shall
denote by $\leq_{\mathrm{c}}$ the componentwise ordering on $\ZZ^m$.

Let $Q$ (resp. $R$) be the set of all linear combinations of the form
$\sum_{j=1}^n\lambda_jg_j$, where all the $\lambda_j$ belong to $\QQ^+$
(resp. $\RR^+$). Hence, $Q\subseteq(\QQ^+)^m$ and
$R\subseteq(\RR^+)^m$.

Note that we also have
$Q=\{x\in\QQ^m\mid(\exists d\in\NN)(dx\in G^+)\}$. Since $G$ is
unperforated, we obtain the following:
 \begin{equation}\label{Eq:QZG+}
 Q\cap\ZZ^m=G^+.
 \end{equation}
Furthermore, the following is an immediate consequence of
Lemma~\ref{L:SysIneq}:
 \begin{equation}\label{Eq:RQQ+}
 R\cap\QQ^m=Q.
 \end{equation}
Now consider elements $a$ and $b$ of $G$ such that $ka\leq b$ (in $G$),
for all $k\in\ZZ^+$. We must prove that
$-a\in G^+$. Note, in particular, that
$b\geq0$ and $b-a\geq0$. Furthermore, by Corollary~\ref{C:Subgrp2},
for all $j\in\fs(n)$, there exists a largest $n_j\in\ZZ^+$ such that
$n_jg_j\leq_{\mathrm{c}}b-a$.

Let $k\in\NN$. Then $(1/k)b-a$ belongs to $Q$, thus there are elements
$\lambda_{k,j}$ of $\QQ^+$ ($j\in\fs(n)$) such that
 \begin{equation}\label{Eq:balambda}
 (1/k)b-a=\sum_{j=1}^n\lambda_{k,j}g_j.
 \end{equation}
Note, in particular, that $0\leq\lambda_{k,j}\leq n_j+1$. Since the
product of the intervals $[0,n_j+1]$, for $j\in\fs(n)$, is a compact
subset of $\RR^n$, there exists an infinite subset $I$ of $\NN$ such
that the sequence
$\vv<\lambda_{k,j}\mid k\in I>$ converges for all $j\in\fs(n)$, say, to
$\lambda_j$. Note that $\lambda_j\geq0$.
Then, taking the limit of the two sides of (\ref{Eq:balambda}) for $k$
going to infinity in $I$ (with respect to the usual topology in
$\RR^m$) yields that $-a=\sum_{j=1}^n\lambda_jg_j$. It follows that
$-a\in R$. However, $-a\in\ZZ^m$, thus, by (\ref{Eq:QZG+}) and
(\ref{Eq:RQQ+}), $-a\in G^+$.
\end{proof}

Compare this result with Corollary~\ref{C:Subgrp2} and
Example~\ref{E:NotArch}.

\begin{remark}
Let $G$ be a finitely presented, unperforated \poag. Suppose, in
addition, that $G$ admits an order-unit, $u$ (this means that $u\geq0$
and $G=\bigcup_{n\in\NN}[-nu,nu]$). By Corollary~\ref{C:FinInt}, the
interval $[-u,u]$ is finite. Therefore, the \emph{order-unit norm}
associated with $u$, defined on $G$ by
 \[
 \|x\|_u=\inf\{p/q\mid p,q\in\NN\text{ and }-pu\leq qx\leq pu\},
 \]
see \cite{Gpoag}, is \emph{discrete}.
\end{remark}

\begin{example}
Endow $G=\ZZ$ with the positive cone $G^+=2\ZZ^+$. Then $G$ is a
finitely presented \poag, and it is $m$-unperforated for every odd
positive integer $m$. However, $G$ is $2$-perforated, thus not
Archimedean.

It is easy to turn this example into a \emph{directed} example
satisfying the properties above, by defining a \poag\ $H$ by
$H=\ZZ\times\ZZ$, and
 \[
 H^+=\{\vv<x,y>\in\ZZ\times\ZZ\mid y=0\Rightarrow x\in 2\ZZ^+\}.
 \]
\end{example}

\section{\eE-ultrasimplicial vector spaces}

The main goal of this section is to prove the following result:

\begin{theorem}\label{T:EUltr}
Let $K$ be a totally ordered division ring, let $E$ be a \povs\ over
$K$. Then every directed finitely presented ordered subspace of
$E$ is contained into a simplicial subspace of $E$.
\end{theorem}

\begin{proof}
Let $\{a_1,\ldots,a_m\}$ be a finite set of generators of $F^+$. By
definition (see Definition~\ref{D:FinPresPoag}), the $K^+$-semimodule
$\SS$ of all elements $\vv<\xi_1,\ldots,\xi_m>\in K^m$ such that
 \begin{equation}\label{Eq:aixii}
 \sum_{i=1}^ma_i\xi_i\geq0
 \end{equation}
is finitely generated. Since $a_i\geq0$
for all $i$, $\SS$ admits a generating subset of the form
 \[
 \{\vv<\lambda_{i,j}\mid 1\leq i\leq m>\mid 1\leq j\leq n\}
 \cup\{\dot e_i\mid 1\leq i\leq m\},
 \]
where $\vv<\dot e_i\mid i\in\fs(m)>$ is the canonical basis of
$K^m$. This can be expressed as follows: the relation
(\ref{Eq:aixii}), with variables
$\xi_1$, \dots, $\xi_m$, is generated by all relations
 \begin{equation}\label{Eq:PresF}
  \begin{aligned}
  \sum_{i=1}^ma_i\lambda_{i,j}\geq0&
  \qquad\text{for all }j\in\fs(n),\\
  a_i\geq0&\qquad\text{for all }i\in\fs(m).
  \end{aligned}
 \end{equation}
By the version for totally ordered division rings of the Grillet,
Effros, Handelman, and Shen Theorem (see Proposition~7.3 of
\cite{Wehr}), there are $p\in\NN$, elements $\alpha_{j,k}\in K$
for $1\leq j\leq n$ and $1\leq k\leq p$, and elements $b_1$, \dots,
$b_p$ of $E^+$ such that
 \begin{equation}\label{Eq:abalpha}
 \begin{pmatrix}a_1\\ \vdots \\ a_m\end{pmatrix}=
 b_1\begin{pmatrix}\alpha_{1,1}\\ \vdots \\ \alpha_{m,1}\end{pmatrix}
 +\cdots+
 b_p\begin{pmatrix}\alpha_{1,p}\\ \vdots \\ \alpha_{m,p}\end{pmatrix},
 \end{equation}
\emph{and} each of the column matrices
$U_k=\begin{pmatrix}\alpha_{1,k}\\ \vdots \\ \alpha_{m,k}\end{pmatrix}$,
for $1\leq k\leq p$, satisfies the system (in $K$) obtained by
substituting in (\ref{Eq:PresF}) the elements of the form
$\alpha_{i,k}$ to the corresponding elements $a_i$. The latter
condition means that the following holds:
 \begin{equation}\label{Eq:PresFun}
  \begin{aligned}
  \sum_{i=1}^m\alpha_{i,k}\lambda_{i,j}\geq0&
  \qquad\text{for all }j\in\fs(n),\\
  \alpha_{i,k}\geq0&\qquad\text{for all }i\in\fs(m),
  \end{aligned}
 \end{equation}
for all $k\in\fs(p)$. Therefore, we have obtained that for all
$\vv<\xi_1,\ldots,\xi_m>\in K^m$ and for all $k\in\fs(p)$,
$\sum_{i=1}^ma_i\xi_i\geq0$ implies that
$\sum_{i=1}^m\alpha_{i,k}\xi_i\geq0$ (it suffices to verify this for
$m$-uples generating the corresponding set of vectors
$\vv<\xi_1,\ldots,\xi_m>$; this holds by (\ref{Eq:PresFun})).
This means that there exists a unique positive homomorphism
$f\colon F\to K^p$ such that the equality
\begin{equation}\label{Eq:faialph}
f(a_i)=\vv<\alpha_{i,1},\ldots,\alpha_{i,p}>
\end{equation}
holds for all $i\in\fs(m)$.
From now on, choose $p$ to be the \emph{smallest} natural integer for
which (\ref{Eq:abalpha}) and (\ref{Eq:PresFun}) are possible, for some
$b_1$, \dots, $b_p$ in $E^+$ and $\alpha_{i,k}$ ($1\leq i\leq m$ and
$1\leq k\leq p$) in $K$. In particular, the following conditions hold.
\begin{enumerate}
\item For all $k\in\fs(p)$, all the elements $\alpha_{i,k}$
($i\leq i\leq m$) belong to $K^+$, and at least one of them belongs to
$K^{++}$ (if $\alpha_{i,k}=0$ for all $i$, then one could decrease $p$
to $p-1$).

\item We have $b_k\in E^{++}$, for all $k\in\fs(p)$ (otherwise, one
could, again, decrease $p$ to $p-1$).

\item No $U_k$, where $1\leq k\leq p$, is a positive linear combination
of the $U_l$, for $l\in\fs(p)\setminus\{k\}$.
\end{enumerate}

Further, by point (i), one can suppose, after an appropriate scaling
of $U_k$, that the following equality
 \begin{equation}\label{Eq:sum=1}
 \sum_{i=1}^m\alpha_{i,k}=1
 \end{equation}
holds, for all $k\in\fs(p)$.
By (\ref{Eq:abalpha}), $F^+$ is contained into the
subspace of $E$ generated by $\{b_1,\ldots,b_p\}$. Thus this also holds
for $F$, because $F$ is directed. Thus, to conclude, it suffices to
prove that $\vv<b_1,\ldots,b_p>$ is the basis of a simplicial subspace
of $E$. Since $b_k>0$ for all $k$, it is easy to see that this is
equivalent to saying that for all $k$, one cannot have
$b_k\propto\sum_{l\ne k}b_l$, where
$x\propto y$ is short for $(\exists\lambda\in K^+)(x\leq y\lambda)$,
for all $x$, $y\in E^+$.

Suppose, to the contrary, that this holds. Without loss of generality,
we can assume that $k=p$, so that
 \begin{equation}\label{Eq:bpbk}
 b_p\propto\sum_{1\leq k<p}b_k.
 \end{equation}
Furthermore, we can, by point (i) above, assume without loss of
generality that $\alpha_{1,p}>0$.

Since $a_1$ is a (positive) linear combination of $b_1$, \dots,
$b_p$, the following is an immediate consequence of (\ref{Eq:bpbk}):
 \begin{equation}\label{Eq:a1bk}
 a_1\propto\sum_{1\leq k<p}b_k.
 \end{equation}
By point (iii) above, no $U_k$ belongs to the
convex hull of the others. By (\ref{Eq:sum=1}), all the $U_k$ belong to
the hyperplane consisting of the column matrices of sum $1$. Therefore,
we can apply Lemma~\ref{L:HyperplHB} to $U_p$ and the convex hull $C$ of
of $\{U_1,\ldots,U_{p-1}\}$. We obtain elements $\beta_1$, \dots,
$\beta_m$ of $K$ satisfying the following properties:
 \begin{align}
 \sum_{i=1}^m\alpha_{i,p}\beta_i=0,&\label{Eq:alpipb}\\
 \sum_{i=1}^m\alpha_{i,k}\beta_i>0,&\qquad
 \text{for all }k\in\fs(p-1).\label{Eq:alkipb}
 \end{align}
Now put $a=\sum_{i=1}^ma_i\beta_i$. By using (\ref{Eq:abalpha}) and
then expanding, we obtain the following:
 \begin{equation*}
 a=b_1\sum_{i=1}^m\alpha_{i,1}\beta_i+\cdots+
 b_{p-1}\sum_{i=1}^m\alpha_{i,p-1}\beta_i+
 b_p\sum_{i=1}^m\alpha_{i,p}\beta_i.
 \end{equation*}
By (\ref{Eq:alpipb}) and (\ref{Eq:alkipb}), every coefficient
$\sum_{i=1}^m\alpha_{i,k}\beta_i$, for $1\leq k\leq p$, is (strictly)
positive, except for $k=p$ where it vanishes. Therefore, we have
obtained that
$\sum_{1\leq k<p}b_k\propto a$. Therefore, by
(\ref{Eq:a1bk}), $a_1\propto a$. Apply to this the homomorphism $f$ of
(\ref{Eq:faialph}). Taking the $p$-th component of the resulting
inequality, we obtain that
 \begin{equation}\label{Eq:alpha1p}
 \alpha_{1,p}\propto\sum_{i=1}^m\alpha_{i,p}\beta_i.
 \end{equation}
However, $\alpha_{1,p}>0$, while, by (\ref{Eq:alpipb}), the right hand
side of (\ref{Eq:alpha1p}) equals $0$; this is a contradiction.
\end{proof}

\begin{example}
The statement of Theorem~\ref{T:EUltr} cannot be extended by removing
the assumption that $F$ is directed. For example, let $E$ be the
\povs\ over $\QQ$ with underlying space $\QQ\times\QQ$, and with
positive cone $E^+=\{\vv<0,0>\}\cup(\QQ^{++}\times\QQ^{++})$ (any
simple, non totally ordered dimension vector space over $\QQ$ would do).
Put $a=\vv<1,-1>$. Then the subspace $F$ of $E$ generated by $a$ is
finitely presented, with $F^+=\{\vv<0,0>\}$.
On the other hand, for any two elements $x$ and $y$ of
$E^{++}$, there exists $n\in\NN$ such that $x\leq ny$ and $y\leq nx$.
Thus every non trivial simplicial subgroup of $E$ is generated by a
single vector of $E^{++}$. Hence, there exists no simplicial subspace
of $E$ containing $F$.
\end{example}

Note the following corollary to Theorem~\ref{T:EUltr}. It gives a
characterization of \eE-ultrasimplicial dimension vector spaces:

\begin{corollary}\label{C:EUltr}
Let $K$ be a totally ordered division ring, and let $E$ be a dimension
right vector space over $K$. Then $E$ is \eE-ultrasimplicial if and only
if it is coherent.
\end{corollary}

The corresponding notions are defined in Definitions \ref{D:UltraS} and
\ref{D:CohOrdGrp}.

\begin{proof}
Suppose first that $E$ is \eE-ultrasimplicial. Let $F$ be a finitely
generated subspace of $E$. By assumption, there exists a simplicial
subspace $S$ of $E$ such that $F\subseteq S$. Since $S$ is simplicial,
it is finitely presented; thus $F$ is also finitely presented (see
Theorem~8.1 of \cite{Wehr}).

Conversely, suppose that $E$ is a coherent dimension right vector space
over $K$. Let $X$ be a finite subset of $E$. Since $E$ is directed, for
all $x\in X$, there are $x_+$ and $x_-$ in $E^+$ such that $x=x_+-x_-$.
Let $F$ be the ordered subspace of $E$ generated by the elements of
the form $x_+$ and $x_-$, for $x\in X$. Note, in particular, that $F$
is a directed subspace of $E$. Since $E$ is coherent, $F$ is finitely
presented. By Theorem~\ref{T:EUltr}, $F$ is contained into a
simplicial subspace $S$ of $E$. Note that $X$ is contained into $S$.
\end{proof}

\begin{corollary}\label{C:EUltrgp}
Let $G$ be a divisible dimension group. Then $G$ is \eE-ultrasimplicial
if and only if it is coherent.
\end{corollary}

The statement that $G$ is \emph{divisible} means, as usual, that
$mG=G$, for all $m\in\NN$.

\begin{proof}
We prove the non trivial direction. So, let $G$ be a coherent,
divisible dimension group. Since $G$ is unperforated and divisible, we
can view $G$ as a \povs\ over $\QQ$. Since $G$ is coherent over $\ZZ$,
it is, \emph{a fortiori}, coherent over $\QQ$. Therefore, by
Corollary~\ref{C:EUltr}, $G$ is \eE-ultrasimplicial as a \povs\ over
$\QQ$. Since $\QQ^n$ is an \eE-ultrasimplicial
dimension group for all $n\in\ZZ^+$ (observe that $\QQ^n$ is
equal to the directed union $\bigcup_{d\in\NN}(1/d)\ZZ^n$),
$G$ is also an \eE-ultrasimplicial dimension group.
\end{proof}

\begin{remark}
Let $K$ be a totally ordered division ring. For a dimension right
vector space $E$ over $K$, the statements that $E$ is an
\eE-ultrasimplicial dimension vector space (resp. dimension group) may
have different meanings. The vector space statement obviously implies
the group statement. The converse is false, as shows the very simple
example of $K=\RR$, and $E=\RR$: then $E$ is an \eE-ultrasimplicial
dimension vector space over $\RR$, but it is not an
\eE-ultrasimplicial dimension group. However, it is easily verified
that if $E$ is a divisible dimension group, then $E$ is
\eE-ultrasimplicial as a dimension group if and only if $E$ is
\eE-ultrasimplicial as a dimension vector space over $\QQ$
(it suffices to observe that all the $\QQ^n$, for $n\in\ZZ^+$,
are \eE-ultrasimplicial as dimension groups).
\end{remark}

\begin{problem}
Is every coherent dimension group \eE-ultrasimplicial?
\end{problem}

Corollary~\ref{C:EUltrgp} provides a positive answer, in the case
of \emph{divisible} \poag s.

\end{document}